\documentclass[11pt,letterpaper]{amsart}
\usepackage{amssymb, amscd, amsfonts, euler, epic,eepic, xypic}
\input xy
\xyoption{all}

\numberwithin{equation}{section}
\def\AA{{\mathbb A}}

\def\CC{{\mathbb C}}

\def\PP{{\mathbb P}}
\def\QQ{{\mathbb Q}} 
\def\RR{{\mathbb R}} 
 
\def\ZZ{{\mathbb Z}}

\def\half{\tfrac{1}{2}}
\def\G{{\Gamma}}
\def\g{{\gamma}}

\def\bb{{\rm bb}}

\def\reg{{\rm reg}}

\def\st{{\rm st}}
\def\ss{{\rm ss}}
\def\sing{{\rm sg}}

\def\bb{{bb}}

\def\bs{\backslash}
\def\bss{{\bs\!\bs}}
\def\pt{{\bullet}}
\def\eps{\epsilon}

\def\Acal{{\mathcal A}}

\def\Ecal{{\mathcal E}} 
\def\Fcal{{\mathcal F}} 
 
\def\Hcal{{\mathcal H}}

\def\Kcal{{\mathcal K}}

\def\Ocal{{\mathcal O}}

\def\Xcal{{\mathcal X}} 
\def\Ycal{{\mathcal Y}} 
\def\Zcal{{\mathcal Z}}

\def\la{\langle}
\def\ra{\rangle}

\newcommand\proj{\operatorname{Proj}}

\newcommand\sym{\operatorname{Sym}}
\newcommand\GL{\operatorname{GL}}

\newcommand\PGL{\operatorname{PGL}}
\newcommand\Res{\operatorname{Res}}
\newcommand\SL{\operatorname{SL}}

\newcommand\supp{\operatorname{supp}}

\newtheorem{theorem}{Theorem}[section]
\newtheorem{lemma}[theorem]{Lemma}
\newtheorem{proposition}[theorem]{Proposition}

\newtheorem{corollary}[theorem]{Corollary}

\theoremstyle{definition}

\theoremstyle{remark} 
\newtheorem{remark}[theorem]{Remark}
\newtheorem{remarks}[theorem]{Remarks}

\begin{document}

\title{The period map for cubic threefolds}
\author{Eduard Looijenga}
\author{Rogier Swierstra}
\email{looijeng@math.uu.nl, swierstr@math.uu.nl}
\address{Mathematisch Instituut\\
Universiteit Utrecht\\ P.O.~Box 80.010, NL-3508 TA Utrecht\\
Nederland}
\thanks{Swierstra is supported by the Netherlands Organisation for Scientific
Research (NWO)}


\begin{abstract} 
Allcock-Carlson-Toledo defined a  period map for cubic threefolds which takes values
in a ball quotient of dimension 10. A theorem of Voisin implies that this is an open embedding. We determine its image and show that on the algebraic level this amounts to identification of the algebra of $\SL (5,\CC)$-invariant polynomials on the representation space $\sym^3(\CC^5)^*$ with an explicitly  decribed algebra of meromorphic automorphic forms on the complex $10$-ball.
\end{abstract}

\maketitle

\section*{Introduction}
The polarized Hodge structure of nonsingular cubic threefold $X\subset \PP^4$ 
is encoded by its intermediate Jacobian, which is a principally polarized 
abelian variety of dimension $5$. It seems hard to characterize the  
abelian varieties that so appear and that is perhaps the reason that 
Allcock-Carlson-Toledo proposed to consider instead the 
cyclic order three cover of $\PP^4$ that is totally ramified over $X$ 
and take the polarized Hodge structure (with $\mu_3$-action) of that. 
This cover  is a smooth cubic fourfold and such a variety has an 
interesting and well-understood Hodge structure: its primitive 
cohomology sits in the middle dimension and has as its nonzero 
Hodge numbers $h^{3,1}=h^{1,3}=1$ and $h^{2,2}=20$. What makes 
this so tractable is that the period map for nonsingular cubic 
fourfolds is very much like that for polarized K3 surfaces: such 
polarized Hodge structures are  parameterized by a bounded symmetric 
domain of type IV in the Cartan classification and of dimension 20 
(so one more than for polarized K3 surfaces). Moreover, the period 
map is a local isomorphism  and even more is true: according to a 
theorem of Voisin, it is an open embedding. Here we are however 
dealing with rather special cubic fourfolds, namely those that 
come with a $\mu_3$-action whose fixed point set is a cubic threefold. 
Their primitive cohomology comes with a $\mu_3$-action as well and one 
finds that such data are parameterized by a symmetric subdomain of the 
type IV  of dimension  10 that is isomorphic to a complex ball. 
Thus the period map of Allcock-Carlson-Toledo is a map from the 
moduli space of  nonsingular cubic  threefolds to a complex $10$-ball 
modulo an arithmetic group. Voisin's theorem implies that this is an 
open embedding and so an issue that remains is the determination of its image.
This is the subject of the present paper. 

Our main result, Theorem \ref{thm:main}, states among other things that 
this image is the complement of a locally symmetric divisor. It also 
identifies  the algebra of $\SL (5,\CC)$-invariant polynomials on the 
representation $\sym^3(\CC^5)^*$ with an algebra of meromorphic automorphic 
forms on the complex $10$-ball. Thus the situation is very much like that 
of the period map for quartic plane curves (which assigns to  such a 
curve the Hodge structure of the $\mu_4$-cover of $\PP^2$ totally ramified 
along that curve). Indeed, this fits in the general framework that we 
developed for that purpose in \cite{ls}.
The same construction applied to a cubic surface yields a cubic threefold with 
$\mu_3$-action (and hence give rise to a cubic fourfold with  
$\mu_3\times\mu_3$-action).  In this manner one can deduce from our 
theorem Allcock-Carlson-Toledo's identification \cite{act} of the 
moduli space of cubic surfaces as a four dimensional ball quotient.
\\

After we finished this paper, we received from Allcock-Carlson-Toledo a manuscript 
\cite{act:ms} in which they also determine the image of the period map. 
Their proof is clearly different from ours.
\\

\emph{Acknowledgements.} It is a pleasure to acknowledge a discussion on 
this topic with Daniel Allcock in the spring of 2001.  We also  
thank Jim Carlson for sending us a draft on the semistable reduction 
of the chordal cubic and its associated fourfold (a preliminary version 
of section 5.3 of \cite{act:ms}).

This paper was written up when the first author was spending the spring term
of 2006 at the Laboratoire J.A.~Dieudonn\'e of the Universit\'e de Nice. He thanks the 
Laboratoire, and especially his host Arnaud Beauville,  for providing such 
pleasant working conditions as well as for  partial support. He is also 
grateful to his own department for giving him leave of absence.
\\

\section{Cubic hypersurfaces of dimension three and four}\label{sect:4folds}

\subsection*{Cubic fourfolds and their residual Hodge lines}
The middle dimensional cohomology group of nonsingular cubic fourfold 
$Y\subset \PP^5$ is  free of rank $23$ and comes with an intersection form that is  
unimodular (because of Poincar\'e duality) and of signature $(21,2)$. If $\eta\in  H^2(Y,\ZZ)$ is the  hyperplane class, then $\la \eta^4,[Y]\ra=3$ and so $3$ is also the self-intersection of $\eta^2\in H^4(Y,\ZZ)$. In particular, the intersection  form is odd.
The primitive part $H_o^4(Y,\ZZ)\subset H^4(Y,\ZZ)$, which  is by definition the orthogonal complement of $\eta^2$, is generated by vanishing cycles. As a 
vanishing cycle has  self-intersection $+2$, this is an even lattice. According to the theory of quadratic forms (see for instance Nikulin \cite{nikulin})  this characterizes the vector $\eta^2$ up to an orthogonal transformation of $ H^4(Y,\ZZ)$ and we  may conclude that we the have a lattice isomorphism 
\[
H_o^4(Y,\ZZ)\cong E_8^2\perp U^2\perp A_2,
\]
where, as usual, $U$ denotes the hyperbolic plane (a lattice spanned by two isotropic vectors which have inner product $1$).

The nonzero Hodge numbers in dimension $4$ are $h^{3,1}(Y)=h^{1,3}(Y)=1$, $h^{2,2}(Y)= 21$ and hence $H^4_o(Y)$ has a Hodge structure of type  IV.
We can represent $H^4_o(Y)$ be means of regular $5$-forms on $\PP^5-Y$. This follows from the fact that $\PP^5-Y$ is affine and the proposition below. 

\begin{proposition}\label{prop:residue}
If $Y$ is a cubic $4$-fold whose singular set is nonsingular of dimension $\le 1$, then we have an exact sequence
\[
0\to H^5(\PP^5-Y)\to H^4(Y_\reg)(-1)\to \ZZ\to 0,
\]
where $H^5(\PP^5-Y)\to H^4(Y_\reg)(-1)$ is the residue map and $H^4(Y_\reg)\to \ZZ$
is integration over a general linear section of dimension two.
\end{proposition}
\begin{proof} 
First consider the exact sequence
\[
H^5(\PP^5)\to H^5(\PP^5-Y)\to H^6_Y(\PP^5)\to H^6(\PP^5)\to H^6(\PP^5-Y).
\]
We have $H^5(\PP^5)=0$ and $H^6(\PP^5-Y)=0$ because $\PP^5-Y$ is affine.
Hence $H^5(\PP^5-Y)$ maps isomorphically to the primitive part of 
$H^6_Y(\PP^5)$. In the exact sequence below 
\[
H^6_{Y_\sing}(\PP^5)\to H^6_Y(\PP^5)\to H^4(Y_\reg)(-1)\to
H^7_{Y_\sing}(\PP^5)
\]
the extremal terms are zero because $Y_\sing$ is smooth of codimension $\ge 4$
and hence the middle map is an isomorphism. The proposition follows.
\end{proof}

Thus an  equation $G\in\CC[Z_0,\dots ,Z_5]$ (a homogeneneous form of degree three) defines an element $[\alpha (G)]\in H^4(Y_\reg,\CC)$  by taking the image of the class 
\[
[\Res_{\PP^5} \frac{dZ_0\wedge\cdots \wedge dZ_5}{G^2}]\in H^5(\PP^5-Y;\CC)
\]
under the map of Proposition \ref{prop:residue}. It will be important for us to verify that
in certain cases $[\alpha (G)]$ is nonzero. In case $Y$ is smooth that is certainly so, for according to Griffiths \cite{griffiths}, $[\alpha (G)]$ is then a generator of $H^{3,1}(Y)$.
It is clear that the span  of $[\alpha(G)]$ in $H^4(Y_\reg;\CC)$ only depends on $Y$.
We write $\Fcal (Y)$ for the one dimensional vector space spanned by 
$G^{-2} dZ_0\wedge\cdots \wedge dZ_5$. We often identify $\Fcal (Y)$  with its image of
under the Griffiths residue  map in $H^4(Y_\reg; \CC)$ if that image is nonzero.

\begin{remark}
On the form level the residue map can be  defined  in terms of the inner product on
$\CC^6$ (or equivalently, in terms of the Fubini-Study metric).
Concretely, consider the normalized  gradient vector field of $G$ in $\CC^6$
\[
N(G):=\frac{\| Z\|^2}{\left\| dG\right\| }\sum_{i=0}^5 \left(\overline{\frac{\partial G}{\partial Z_i}}\right)\frac{\partial }{\partial Z_i}.
\]
Letting $\iota$ stand for the inner product (the contraction of a vector with a form),
then \[
\iota_{N(G)} d\iota_{N(G)} \left(\frac{dZ_0\wedge\cdots \wedge dZ_5}{G^2}\right)
\] 
is a $5$-form on $\CC^6-\{ 0\}$ of Hodge level $\ge 4$ on 
whose restriction to the zero set of $G$ is closed. It is also invariant under scalar multiplication and so it has a residue 
at infinity: this is a closed $4$-form on $\PP^5$ whose restriction $\alpha(G)$ to $Y_\reg$ is closed.  It is a sum of  a form of type $(3,1)$ and one of type $(4,0)$ that represents the  above residue up to a universal nonzero scalar. We might therefore also think of 
$\Fcal (Y)$ as a line of forms on $Y_\reg$ (which however depends on the Fubini-Study metric).
\end{remark}

The  simple hypersurface singularities in dimension four, $A_k$, $D_{k\ge 4}$, $E_6$, $E_7$, $E_8$, are the `double suspensions' of the Kleinian singularities that bear 
the same name. We recall that a suspension of a hypersurface singularity adds to its 
equation a square in a new variable. Doing this twice does not affect the 
monodromy group of its miniversal deformations and so in the present case
we have finite monodromy groups: such singularities go largely unnoticed 
for the period map. 

We wish to establish that for some reduced cubic fourfolds $Y$,
$(Y,\Fcal (Y))$ is a \emph{boundary pair} (in the sense of \cite{ls}) and we will also want to know its type. Here we use that notion in a slightly more general  sense than Definition (2.3) of \emph{op.\ cit.} in that we  allow (at least in principle) the possibility that $\Fcal(Y)$ maps to zero in $H^4(Y_\reg;\CC)$, but then require that the map from 
$H^o_4(Y_\reg;\CC)$ (or the part that matters to us---in the present case an eigenspace for an action of the cyclic group of order three)  to the primitive homology of a smoothing of $Y$  be nontrivial. In that case $Y$ still imposes via Lemma 1.2 of \cite{ls} a nontrivial linear constraint on the limiting behavior of the period map. 

The first class of examples is furnished by the isolated hypersurface  
singularities in dimension four that come after the simple ones:  the double suspensions of the simple-elliptic singularities 
$\tilde E_6$, $\tilde E_7$, $\tilde E_8$.

\begin{proposition}\label{prop:type2}
Let $Y$ be a  cubic fourfold with a singular point of type $\tilde E_6$, $\tilde E_7$ or  
$\tilde E_8$. Then $(Y,\Fcal(Y))$ defines a  boundary pair of type II.
\end{proposition}
\begin{proof}
It is enough to verify that if  $G$ is an equation for $Y$, then $\alpha (Y)\wedge \bar\alpha (Y)$ is not integrable and that $H_4(Y_\reg)$ contains 
an isotropic lattice of rank two that is mapped by $\alpha (Y)$ to a lattice in $\CC$. We will see that this is essentially a local issue 
that has been dealt with in singularity theory. Let us for concreteness 
assume that $Y$ has a singularity $o$ of type $\tilde E_8$. We may choose 
local complex-analytic coordinates 
$(z_1,\dots ,z_5)$ at $o$ such that $Y$ is given there as 
$f(z)=z_1^6+z_2^3+\lambda z_1z_2z_3+z_3^2+z_4z_5$. 
Notice that $f$ is weighted homogenenous of degree $6$ with  weights 
$(1,2,3,3,3)$.  The  residue $\alpha$ of $f^{-2}dz_1\wedge\cdots \wedge dz_5$ on the smooth part of the zero set of $f=0$ is homogeneous of degree zero
(equivalently: $\CC^\times$-invariant). The form $\alpha\wedge\bar\alpha$ is positive everywhere and $\CC^\times$-invariant also, and hence will not be  integrable near $o$. It is well-known (and implied by the work of Steenbrink \cite{steenbrink}) that  the link $L$ of this singularity has the property that $H_4(L)$ is free of rank two and that the periods of $\alpha$ on it are the periods of 
of the elliptic curve defined by $z_1^6+z_2^3+\lambda z_1z_2z_3+z_3^2$ in 
a weighted projective space.  It is also known (see \cite{steenbrink}) that 
if we multiply $\alpha$ with an element of the maximal ideal 
of $\CC\{z_1,\dots,z_5\}$, then it becomes exact on the germ of $Y_\reg$ at 
$o$. Since  $\alpha (Y)$ equals $\alpha$ up to a unit in $\CC\{z_1,\dots,z_5\}$, the image of 
$H_4(L)\to H_4(Y_\reg)$ is as required and the proposition follows.
\end{proof}

In the following lemma we use the notion of boundary pair in the above more general sense.

\begin{lemma}\label{lemma:type1}
Let $Y$ be a cubic fourfold whose singular locus has as an irreducible component a curve such that $Y$ has a transversal singularity  of type $A_2$ along the generic point of that curve. If the primitive homology $H^0_4(Y_\reg)$ has nontrivial intersection pairing, then $(Y,\Fcal (Y))$ defines a boundary pair of type I and $H^0_4(Y_\reg)$ is positive semi-definite.
\end{lemma}
\begin{proof} 
Choose complex-analytic coordinates $(z_1,\dots ,z_5)$ at a generic point of the curve in question such that $Y$ is there given by $f(z)=z_1^3+z_2^2+z_3^2+z_4^2$. So $f$ is weighted homogenenous of degree $6$ with  weights 
$(2,3,3,3)$.  The  residue $\alpha$ of $f^{-2}dz_1\wedge\cdots \wedge dz_5$ on the smooth part of the zero set of $f=0$ is homogeneous of degree $-1$ and 
hence the form $\alpha\wedge\bar\alpha$  will not be integrable near the origin. Argueing as in the proof of Proposition \ref{prop:type2} we find that $\alpha (G)$ is not integrable. 
Now let $Z$ be a $4$-cycle on $Y_\reg$ that is  perpendicular to the hyperplane class and has nonzero self-intersection number $Z\cdot Z$. If $\Ycal/\Delta\subset \PP^5_\Delta$ is any smoothing  of $Y$, then $Z$ extends as a relative cycle $\Zcal/\Delta$. Clearly,  $\int_{Z_t} \alpha(Y_t)$ is bounded. On the other hand, 
$\int_{Y_t}\alpha (Y_t)\wedge \overline{\alpha(Y_t)}$ tends to infinity as $t\to 0$.
This implies that any limiting point of the line in $H^4(Y_t;\CC)$ spanned by 
$\alpha(Y_t)$ is in the hyperplane defined by $[Z_t]$. Since $Z_t\cdot Z_t\not=0$,
this hyperplane must be of type I: $Z_t\cdot Z_t>0$.
\end{proof}

\begin{remark}
In the preceding lemma the assumption that $Y_\sing$ contains a curve along which we have a transversal $A_2$-singularity can be weakened to: $Y_\sing$ contains an irreducible component of dimension one (the above argument also works for transversal singularity type $A_1$ and hence for any singularity type that is worse).
\end{remark}

The GIT of cubic fourfolds is probably not sufficiently worked out yet to make 
it feasible at present to verify whether every semistable cubic fourfold
yields a boundary pair, but
we shall see that we can do this for cubic fourfolds attached to semistable 
cubic threefolds.

Let $X\subset\PP^4$ be a nonsingular cubic threefold. Following Allcock-Carlson-Toledo its period map is best studied by passing to the $\mu_3$-cover $Y\to\PP^4$ which ramifies over $X$. This cover is a cubic $4$-fold. To be more precise, let 
an equation for $X$ be $F\in\CC[Z_0,Z_1,Z_2,Z_3,Z_4]$. Then $G:=F-Z_5^3$ is an equation for  $Y$.  Moreover, $H^{3,1}(Y)$ comes with the generator $[\alpha (G)]$.

The GIT for cubic hypersurfaces $X\subset\PP^4$ has been carried out 
independently by Allcock \cite{allcock} and Yokoyama \cite{yokoyama1}. 
They find that such a $X\subset \PP^4$ is stable if and only its 
singularities are of type $A_1$, $A_2$, $A_3$ or $A_4$. This means that 
$Y$ has singularities of type $A_2$, $D_4$, $E_6$ or $E_8$ respectively 
(add a cube in a new variable).  The minimal strictly semistable cubic threefolds $X$ are the following:
\begin{enumerate}
\item[($D_4^3$)] $X_\sing$ consists of three $D_4$-singularities. Such an $X$ is unique up a linear transformation. The associated fourfold $Y$ has three $\tilde E_6$ singularities (double suspensions of degree three simply-elliptic singularities of  zero $j$-invariant). 
Hence $Y$ yields a boundary pair of type II.
\item[($A_5^2$)] $X_\sing$ consists of two $A_5$-singularities, perhaps  
also with a singularity of type $A_1$. This makes up a one  parameter family of 
$\PGL (5)$-orbits. The associated fourfold $Y$ has two $\tilde E_8$ singularities 
(in fact double suspensions of degree one simply-elliptic singularities 
of  zero $j$-invariant), and possibly an $A_2$-singularity. Hence $Y$ yields a 
boundary pair of type II.
\item[($A_1^\infty$)] $X$ is a chordal cubic: it is the secant variety 
of a rational normal curve in $\PP^4$ of degree $4$; this curve equals 
$X_\sing$ and the transversal singularity is of type $A_1$. It lies in the the closure of
the curve that parameterizes the  $(A_5^2)$-case.
\end{enumerate}
So the GIT boundary consists of  an isolated point ($D_4^3$) and an irreducible curve
that is the union of ($A_5^2$) and ($A_1^\infty$).

\section{Eisenstein lattices}

\subsection*{Generalities}
We fix a generator $T$ of $\mu_3$ so that the group ring $\ZZ\mu_3$ is identified with
$\ZZ[T]/(T^3-1)$. This identifies the number field $\QQ(\mu_3)$ with 
$\QQ [T]/(T^2+T+1)$ and its ring of integers with 
$\ZZ[T]/(T^2+T+1)$ (which is therefore a quotient of $\ZZ\mu_3$). The latter is called  the \emph{Eisenstein ring} and we shall denote it by $\Ecal$.
If we substitute for $T$ the standard choice (relative to a choice of $\sqrt{-1}$) of a primitive 3rd root of unity, $\zeta=-\half+\half\sqrt{-3}$ , then $\Ecal$ gets identified with the set of $\half (a+b\sqrt{-3})$ with $a,b\in\ZZ$ of the same parity.

If $\mu_3$ operates on a finitely generated free abelian group $A$, then $\Ecal\otimes_{\ZZ\mu_3}A$ can be identified with the quotient of $A$ by the fixed point subgroup
$A^{\mu_3}$. And if the latter happens to be trivial (so that $A$ is a $\Ecal$-module), then 
$\CC\otimes A$ splits according to  the characters of $\mu_3$ as
\[
\CC\otimes A=(\CC\otimes A)_\chi \oplus (\CC\otimes A)_{\bar{\chi}},
\]
where $\chi:\mu_3\subset\CC^\times$ is the  tautological character. The first summand
may be identified with $\CC\otimes_\Ecal A$ and the second summand is the complex conjugate of the first. If $A$ also comes
with an integral  $\mu_3$-invariant symmetric bilinear form $(\, \cdot \,) :A\times A\to\ZZ$, then 
\[
\phi : A\times A\to\Ecal, \quad \phi (a,a')=-(a\cdot a')\zeta +(a\cdot Ta')
\]
is  skew-Hermitian over $\Ecal$. It is such that $\phi (a,a)=-\half\sqrt{-3} (a\cdot a)$
(so $(\, \cdot \,)$ had to be even).  Multiplication by $\sqrt{-3}$ turns $\phi$ 
into a Hermitian form
\[
h(a,a'):=\sqrt{-3}\phi (a,a')=\tfrac{3}{2}(a\cdot a') +  
\sqrt{-3} (a\cdot \half a'+ T a')
\]
with $h(a,a)=\frac{3}{2}(a\cdot a)$. Conversely, every finitely generated torsion free $\Ecal$-module equipped with an 
$\Ecal\sqrt{-3}$-valued Hermitian form  (or equivalently, an $\Ecal$-valued skew-Hermitian form $\phi$) so arises and that is why we call these data an \emph{Eisenstein lattice}.

We shall be concerned with certain  Eisenstein lattices denoted $\Lambda_k$ and so we recall their definition: $\Lambda_k$ is a free $\Ecal$-module  with generators $r_1,\dots ,r_k$, whose  Hermitian form is characterized by
\begin{equation}
h(r_i, r_j)=
\begin{cases}
3 &\text{ if } j=i,\\
\sqrt{-3} &\text{ if } j=i+1,\\
0 &\text{ if } j>i+1.
\end{cases}
\end{equation}
This is equivalent to $r_i\cdot r_i=2$, $r_i\cdot r_{i+1}=0$,  $r_i\cdot Tr_{i+1}=1$
and  $r_i\cdot T^kr_j=0$ for $j>i+1$ and all $k$.
This lattice is isomorphic to its conjugate, for  the matrix of $h$ on the basis 
$((-1)^ir_i)_i$ is conjugate to the matrix of $h$ on $(r_i)_i$.

Here are a few cases of special interest to us. The $A_2$-lattice has just two rotations of order three which are each others inverse (these are also its Coxeter transformations) and a $\mu_3$-action thus obtained  it becomes an Eisenstein lattice 
isomorphic to $\Lambda_1$. If we do something similar to an $E_8$-lattice by letting
$T\in\mu_3$ act as the  tenth power of a Coxeter transformation
(which has order $30$), then the resulting Eisenstein lattice is isomorphic to 
$\Lambda_4$. The hyperbolic Eisenstein lattice $U_\Ecal$ is by definition spanned by two isotropic vectors with inner product $\sqrt{-3}$; its underlying integral lattice is 
$U^2$. It is known that $\Lambda_{10}\cong \Lambda_4\perp U_\Ecal\perp 
\Lambda_4$. 

\subsection*{The vanishing Eisenstein lattice attached to a chordal cubic} 
If $X\subset \PP^4$ is a nonsingular cubic $3$-fold and $\PP^5\supset Y\to \PP^4$ the associated cubic $4$-fold with $\mu_3$-action,  then the $\mu_3$-invariant part of $H^4(Y;\QQ)$ can be identified with $H^4(\PP^4;\QQ)$. This is also the image of 
$H^4(\PP^5;\QQ)\to H^4(Y;\QQ)$ and hence is spanned by $\eta^2$. It follows that
$H^4_o (Y,\ZZ)$ is in a natural manner a $\Ecal$-module. The  intersection pairing turns it into an Eisenstein lattice. We  use the degeneration of $X$  into a chordal cubic to determine its isomorphism type.

Jim Carlson has determined the limiting Hodge structure 
for a linear smoothing of the chordal cubic as well as for the associated smoothing
of cubic fourfolds. The lemma below can also be derived from his computations.

\begin{lemma}\label{lemma:vanlattice}
Let $X\subset \PP^4$ be the chordal cubic, $\Xcal/\Delta\subset \PP^4_\Delta$ a general linear smoothing of $X$ over the unit disk $\Delta$ and $X'$ a general fiber of this smoothing. Denote by 
$Y\subset \PP^5$, $\Ycal/\Delta\subset\PP^5_\Delta$  and $Y'$ 
the associated (relative) cubic fourfolds. Then the kernel of a natural map $H_4(Y';\Ecal)\to H_4(Y;\Ecal)$ equipped with the intersection pairing and the $\mu_3$-action contains an Eisenstein lattice isomorphic to $\Lambda_{10}$. 
\end{lemma}
\begin{proof}
By assumption the smoothing of $X$ has the form $F_t=F+tF'$, where $F$ is an equation for $X$. We suppose that  that $(F'=0)$ meets  the singular set $C$ of $X$ transversally and that for $0<|t|<1$, $X_t=(F_t=0)$ is smooth. For such $t$, $C\cap X_t$ consists of $4.3=12$ distinct points. Near a point of $C-C\cap X_t$ resp.\ $C\cap X_t$ we can find local analytic coordinates 
$(z_1,\dots ,z_5)$ on $\PP^5$ such that the smoothing $\Ycal_\Delta$ is given by $z_1^3+z_2^2+z_3^2+z_4^2=t$ resp.\ $z_1^3+z_2^2+z_3^2+z_4^2=tz_5$ with $\mu_3$ affecting only the first coordinate.

Choose an oriented embedded circle  $\gamma$ in $C$ that contains $C\cap X_t$, label the points of $C\cap X_t$ in a corresponding cyclic manner $\{ p_i\}_{i\in\ZZ/12}$, and denote by $\gamma_i$ the part of $\g$ that goes from $p_{i-1}$ to $p_i$. In what follows, only the isotopy class of $\g_i$ matters.

For $\eps>0$ small and $i\in\ZZ/12$ given, we construct over $\gamma_i$ a cycle 
$\G_i$ in $Y':=Y_\eps$ as follows. Since this essentially only involves the  topology
we may suppose that  $\g_i$ is very small so that $p_{i-1}$ and $p_i$  are about to coalesce. This allows us to assume  that 
$\gamma_i$ is contained in a complex-analytic coordinate patch $(U; z_1,\dots ,z_5)$
such that $C\cap U$ is open in the $z_5$-axis, $z_5(p_{i-1})=-1, z_5(p_i)=1$
and the smoothing is of the form 
\[
z_1^3+z_2^2+z_3^2+z_4^2=t(1-z_5^2),
\]
with as before, $\mu_3$ affecting  the first coordinate only,  and that $\g |U$ is simply the intersection of $U$ with the real part of the $z_5$-axis (so that $\g_i$ is the interval
$[-1,1]$ in that axis). 

Let $D_i$ be the chain on $Y'$ defined by 
\[
x_1^3+x_2^2+x_3^2+x_4^2=\eps(1-x_5^2), \quad x_i\in\RR,\; x_1\ge 0,\;  -1\le x_5\le 1.
\]
We notice that $D_i$ is a topological $4$-disk, for its projection on the $x_5$-axis 
has as fiber over $x_5\in [-1,1]$ the $3$-disk in case $|x_5|\not=1$
(with $\eps(1-x_5)^2-x_1^3$ as the radial parameter and 
$(x_2,x_3,x_4)/|(x_2,x_3,x_4)|$ as angular parameter), 
and as fiber over $\pm1$ the singleton  $\{(0,0,0,0, \pm 1)\}$. 
We orient $D_i$ by taking the orientation defined by the coordinates ($x_1,\dots ,x_4)$ at  the point $(\sqrt[3]{\eps},0,0,0,0)$ so that we may regard it as a chain. 
It has the same boundary as $T D_i$ and hence $\G_i:=(1-T) D_i$ is the cycle defined
by an oriented $4$-sphere. 
\\

\emph{Claim 1.  The intersection numbers $(\G_i\cdot T^k\G_{i+1})_{k\in\ZZ/3}$
are up to a cyclic permutation equal to $(1,0,-1)$.}

Near $(0,0,0,0,1)$, $(z_1,\dots ,z_4)$ is a coordinate system for $Y'$. 
In terms of that system, $D_i$ is simply the set for which all coordinates are real and
$x_1\ge 0$. Hence $\G_i$ is there defined by $z_2,z_3,z_4$ real and $z_1$ or 
$\zeta^{-1}z_1$ real and $\ge 0$. After a possible isotopy, the sphere $D_{i+1}$ meets
$U$ in the set where where $z_1$ is a primitive $6$th root of unity and $z_2,z_3,z_4$ 
are purely imaginary (so that $1-z_5^2$ is real and $\le 0$). These are oriented topological submanifolds  meeting in $(0,0,0,0,1)$ only. Their intersection
number is the same as that of an arc $a$ in $\CC$ composed of two half rays making 
an angle $2\pi/3$ and the transform of $-a$ under a primitive $6$th root of unity. 
It is clear that if we replace $\G_{i+1}$ by $T^k\G_{i+1}$, then we must multiply 
this primitive $6$th root of unity by $\zeta_3^k$. Now the claim follows from 
the easily seen fact that if $\zeta_6=\half +\half\sqrt{-3}$ and $a'$ is $\zeta_6 a$ 
with its opposite orientation, then 
$a\cdot a'=1$, $a\cdot \zeta_3 a'=0$ and $a\cdot\zeta^2_3 a'=-1$.
\\

\emph{Claim 2.  $\G_i$ generates in $H_4(Y')$ a copy of $\Lambda_1$: 
$(1+T+T^2)\G_i=0$ and $\G_i\cdot \G_i=2$.}

It is clear that $(1+T+T^2)\G_i=0$.

On $U\cap Y'$  we have a flow of the form 
\[
\Phi_\alpha (z_1,\dots ,z_5)= 
(e^{\sqrt{-1}\alpha/3}z_1,e^{\sqrt{-1}\alpha/2}z_2,e^{\sqrt{-1}\alpha/2}z_3,
e^{\sqrt{-1}\alpha/2}z_4, \phi_\alpha(z_5)),
\]
where $\alpha$ is small and $1-\phi_\alpha (z_5)^2=e^{\sqrt{-1}\alpha}(1-z_5^2)$ (since $1-z_5^2$ 
has only simple zeroes this is well defined). It has $(0,0,0,0,\pm 1)$ as fixed points.
We see that for nonzero $\alpha$, $\Phi_\alpha(\G_i)$ meets $\G_i$ in these fixed points only.
The argument above shows that at each of these points the intersection number is 1:
it is the intersection number of the real part of  $\CC^4$ with its transform under 
$\Phi_\alpha$ (as acting on the first four coordinates).
\\

It is clear that if $j-i\not=\pm 1$, then $\G_i\cdot T^k\G_j=0$  for every $k$. 
Upon replacing $\G_{i+1}$ for  $i=1,\dots ,9$ successively by an element of $\pm\mu_3 \G_{i+1}$, we can arrange that  $\G_i\cdot \G_{i+1}= 0$ and $\G_i\cdot T\,\G_{i+1}= 1$ for $i=1,\dots , 9$ so that $\G_1,\dots ,\G_{10}$ span in $H_4(Y';\Ecal)$ a sublattice of type  $\Lambda_{10}$.
 \end{proof}

The following proposition plays a central role in this paper and will be proved in Section \ref{sect:chordal}.

\begin{proposition}\label{prop:chordal}
Let $X\subset\PP^4$ be a  chordal cubic and $Y$ the cubic $4$-fold that is a $\mu_3$-cover ramified over $X$. Then $H_\pt (Y_\reg ;\CC)_\chi=H_4 (Y_\reg ;\CC)_\chi$ and the latter is of dimension one. Moreover, the intersection pairing defines a positive Hermitian form on  this space so that $(Y,\Fcal(Y))$ is a boundary pair of type I.
\end{proposition}

\begin{remark}
Let $F\in\CC[Z_0,\dots ,Z_4]$ define  the chordal cubic and put $G:=F-Z_5^3$. If there is a simple way to prove that the $6$-form  $G^{-2}dZ_0\wedge\cdots \wedge dZ_5$ is not exact, then  Proposition \ref{prop:residue}  implies that the Griffiths residue $\alpha (G)$ defines
a nonzero class in $H^4(Y_\reg;\CC)$. We then could use that fact to 
bypass the preceding proposition, resulting in the elimination of Section \ref{sect:chordal}
and hence in a substantial shortening of proof of our main Theorem \ref{thm:main}.
\end{remark}

Here is an interesting corollary to Lemma \ref{lemma:vanlattice}.
It seems to suggest that $H_4(Y_\reg;\Ecal)$ is isomorphic to $\Lambda_1$.

\begin{corollary}\label{cor:elattice}
For $Y'$ as in  Lemma \ref{lemma:vanlattice},  $H_o^4(Y')$ is an Eisenstein lattice
and is as such isomorphic to $\Lambda_{10}\perp\Lambda_1$.
\end{corollary}
\begin{proof}
Since $H^4(Y';\CC)^{\mu_3}\cong H^4(Y'{}^{\mu_3};\CC)=H^4(\PP^4;\CC)$ is spanned by $\eta^2$, it follows that $H_o^4(Y')$ is an Eisenstein lattice. 
The integral lattice underlying $\Lambda_{10}$ is $E_8^2\perp U^2$, hence is unimodular. If we identify $H^4(Y')$ with $H_4(Y')$ via Poincar\'e duality, and use the fact that $\Lambda_{10}$ is isomorphic to its conjugate, 
then the lemma in question proves  that $H^4(Y')$ contains a copy of  $\Lambda_{10}$.

The orthogonal complement of this copy of $\Lambda_{10}$ in $H^4(Y')$
is unimodular, odd, positive definite and  of rank $3$ and hence isomorphic to $\ZZ^3$ equipped with its standard form $x_1^2+x_2^2+x_3^2$. 
The vectors of with self-product $3$ are those that have  each coordinate $\pm 1$, hence lie all in the same orbit of the integral orthogonal group.
We may therefore arrange that the isomorphism takes $(1,1,1)$ to $\eta^2$, so that the orthogonal complement of the $\Lambda_{10}$-copy in $H^4_o(Y')$ is identified with the orthogonal complement of $(1,1,1)$ in $\ZZ^3$. The latter is an $A_2$-lattice and any $\Ecal$-structure on that lattice makes it isomorphic to $\Lambda_1$.
The corollary follows.
 \end{proof}

\section{The main result}

In this section we fix a complex vector space $U$ of dimension $5$ and 
abbreviate the $\GL (U)$-representation  $\sym^3U^*$  by $S$.
Let $\Xcal\subset \PP(U)_S$ the universal cubic. The latter is given by a single equation $F\in \CC[S\times U]$. Denote by  
$f:\Ycal\subset \PP_S(U\oplus \CC)\to S$ the $\mu_3$-cover defined by
$G:=F+w^3$. It is invariant under the obvious $\SL (U)$-action. If we fix a generator $\mu\in \det (U^*)$, then the Griffiths residue
construction applied to $G^{-2}\mu$ yields a relative $(3,1)$-form $\alpha$
on the part $\Ycal^\circ$ of $\Ycal$ where $\Ycal$ is smooth over $S$. 
We denote by $S^\circ$ the locus where $\Ycal$ is smooth over $S$ so that $\Ycal_{S^\circ}\subset\Ycal^\circ$.
\\

\subsection*{The refined period map} 
We fix an odd unimodular  lattice $L$ of signature $(21,2)$, a vector $v_o\in L$ with $v_o\cdot v_o=3$ whose orthogonal complement  $L_o$ is even, and a $\mu_3$-action on $L$ whose fixed point set is $\ZZ v_o$ and for which $L_o$ isomorphic to 
$\Lambda_{10}\perp\Lambda_1$ as a $\Ecal$-lattice. We shall write $\Lambda$
for $L_o$ as an $\Ecal$-lattice. The quadratic form on $L_o$ becomes  a Hermitian form on $\Lambda$ that takes values in $\sqrt{-3}\Ecal$; we denote that form by $h$. We denote the underlying $\CC$-vector space of $\Lambda$ by
$H$ (so $H:=\CC\otimes_\Ecal  \Lambda$ and $\CC\otimes L_o=H\perp \bar H$) and we let $H_+\subset H$ be the set 
of $v\in H$ with $h (v,v)<0$. Then $\PP(H_+)$ is the symmetric domain of
the unitary group of $H$ and is isomorphic to a complex  $10$-ball. The restriction of 
$\Ocal_{\PP(H)} (-1)$ to $\PP(H_+)$  is our basic automorphic line bundle (its $11$th tensor power is the equivariant canonical bundle of $\PP(H_+)$);  we will 
denote that line bundle by $\Acal(1)$. The subgroup $\Gamma\subset O(L)$ of $\mu_3$-automorphisms of $L$ that fix $v_o$ is arithmetic  and acts properly on
$H_+$ and $\PP(H_+)$. The Baily-Borel theory asserts that 
\[
\oplus_{k\ge 0} H^0(\PP(H_+),\Acal (k)))^\Gamma
\]
is a finitely generated graded algebra (of automorphic forms) whose $\proj$ defines a normal projective completion $\G\bs \PP(H_+)\subset \G\bs \PP(H_+)^\bb$  of the orbit space by adding one point (a \emph{cusp}) for every $\G$-orbit in $\partial\PP(H_+)$ whose points are defined over $\QQ(\zeta)$. 
\\

We shall interpret $\PP(H_+)$ as a classifying space for certain  Hodge structures on 
$L_o$ of weight $4$ polarized by $(\, \cdot \, )$ and invariant under $\mu_3$:
giving a  point of $\PP(H_+)$ amounts to giving a line $F^3\subset H$ on which $h$ is negative and such a line determines a weight four polarized Hodge structure on $L_o$:
\[
\CC\otimes_\ZZ L_o=H^{3,1}\perp H^{2,2}\perp H^{1,3},
\]
with $H^{3,1}:=F^3$, $H^{1,3}:=\overline{F^3}$ and  $H^{2,2}$ the orthogonal complement of $F^3\perp\overline{F^3}$. Conversely, any Hodge structure on 
$L_o$ with  $(h^{3,1}, h^{2,2},h^{1,3})=(1,20,1)$, 
polarized by  the quadratic form and with $H^{3,1}$ in the eigen space of the 
tautological character $\chi:\mu_3\subset\CC^\times$ is thus obtained.  
 
We have of course arranged that if $Y\subset \PP (U\oplus \CC)$ is a smooth fiber of
$f$, then  there is a $\mu_3$-isomorphism $H^4(Y)\cong L$ that takes $\eta^2$ to $v_o$ (this follows from Corollary \ref{cor:elattice}). Such  isomorphisms (also called \emph{markings} of $Y$) are permuted simply transitively by  $\Gamma\subset O(L)$.
A marking carries $\alpha (Y)\in H^{3,1}(Y)$ to a point of $H_+$. This defines a refined period  map $ S^\circ\to \Gamma\bs H_+$.
It is evidently constant on the $\SL (U)$-orbits. But there is in fact $\GL (U)$-equivariance: if $F\in S^\circ$, then let $F'$ be its transform under the scalar $t\in \CC^\times\subset \GL (U)$: $F'=t^{-3}F$. If  $(z',w')=(tz,w)$, then clearly,
$F'(z')-(w')^3=F(z)-w^3$ and hence we get an isomorphism $Y_F\cong Y_{F'}$.
This isomorphism pulls back $(F'-(w')^3)^{-2}\mu$ to 
$t^{5}(F-w^3)^{-2}\mu$ and hence sends $\alpha (F')$ to $t^{5}\alpha (F)$.
Thus the refined period map defines a  morphism
\[
P: \GL (U)\bs S^\circ \to  \Gamma\bs \PP(H_+)
\]
that is covered by morphism of line bundles that sends $p^*\Acal $ to 
the fifth power of the line bundle over the left hand side defined by the determinant character
$\det: \GL(U)\to\CC^\times$ (which is the geometric quotient of $\Ocal_{S^\circ}(3)$ by $\SL (U)$). 
In this way we get a $\CC$-algebra homomorphism from 
$\oplus_{k\ge 0} H^0(\PP(H_+),\Acal (3k))^\Gamma$ to 
the part of $\CC [S]^{\SL(U)}$ spanned by the summands of degree a multiple of $5$. Since the center of $\SL(U)$ is $\mu_5$ and acts on $\sym^3 U^*$ faithfully by scalar multiplication,  $\CC [S]^{\SL(U)}$ only  lives  in degrees that are multiples of $5$.
For a similar reason, $\oplus_{k\ge 0} H^0(\PP(H_+),\Acal (k))^\Gamma$ only lives in degrees that are multiples of $3$: the center of $\G$ contains $\mu_3$, which acts in the obvious manner on $\Acal$ as scalar multiplication. Thus we find a $\CC$-algebra homomorphism 
\[
p: \oplus_{k\ge 0} H^0(\PP(H_+),\Acal (k))^\Gamma\to \CC [\sym^3 U^*]^{\SL(U)}
\]
that multiplies degree by $\frac{5}{3}$. The $\proj$ of the right hand side is 
the GIT compactification of $\GL (U)\bs S^\circ$, namely $\GL(U)\bss S^\ss$, where
$S^\ss$ denotes the semistable locus and $\bss$ indicates that we take the geometric quotient. It contains  $\GL (U)\bs S^\st$ as a dense-open subset,  where $S^\st$ denotes the $\SL(U)$-stable locus in $S$. Following Allcock \cite{allcock} $S^\st$ is precisely the set of $s\in S$ for which the fiber $Y_s$ has only simple singularities in the sense of Arnol'd (double suspensions of Kleinian singularities). These have the property that the of monodromy of the fibration 
 $\Ycal_{S^\circ}/S^\circ$ near such a fiber is finite. It is well-known \cite{gt} that this 
implies that the period mapping extends across such singularities as a map
\[
P: \GL (U)\bs S^\st \to  \Gamma\bs \PP(H_+).
\]
A local Torelli theorem tells us that $P$ is a local isomorphism. 
The GIT boundary of $\GL (U)\bs S^\st $ is of dimension one and has as a distinguished
point  the isomorphism type of the chordal cubic. Away from
this point the variation of polarized Hodge structure defined by our family  has degenerations of  type II only.
 
\subsection*{The image of the period map} We recall from  \cite{ls} that 
a boundary pair defines a $\G$-orbit $\Kcal$ of hyperplanes of  $H$. If it is of type I, then each  
$K\in \Kcal$ meets $H_+$. Otherwise it is of type II and $\PP(K)$ intersects 
the closure of $\PP(H^+)$ in a single point only; this point lies on the boundary and is defined over $\QQ(\zeta)$ (hence defines a cusp).

\begin{theorem}\label{thm:main}
The period map defines an isomorphism from the moduli space of stable cubic $3$-folds onto a $10$-dimensional ball quotient minus an irreducible locally symmetric
hypersurface. To be precise, let $\Hcal$ be the collection of hyperplanes in $H$
that are the complex span of an Eisenstein sublattice of $\Lambda$ isomorphic to 
$\Lambda_{10}$. Then $\Hcal$ is a $\G$-orbit and
if we write $H^\circ_+$ for $H_+-\cup_{K\in\Hcal} K_+$, then 
$P$ maps $GL (U)\bs S^\st$ isomorphically onto $\G\bs \PP(H^\circ_+)$. 
Moreover, $P$  induces an isomorphism (of degree $\frac{5}{3}$) from the $\CC$-algebra of meromorphic $\G$-automorphic forms 
$\oplus_{k\in \ZZ} H^0(H^\circ_+, \Acal (k))^\G$ (we allow arbitrary poles along
the hyperplane sections indexed by $\Hcal$) onto the $\CC$-algebra of invariants
$\CC [\sym^3 U^*]^{\SL (U)}$. In particular, $H^0(H^\circ_+, \Acal (k))^\G=0$ for
$k<0$ and the GIT compactification $\GL (U)\bss S^\ss$ gets identified with 
$\proj \oplus_{k\ge 0} H^0(H^\circ_+, \Acal (k))^\G$.
\end{theorem}
\begin{proof}
Proposition \ref{prop:chordal} shows that the chordal cubic defines a boundary pair 
$(Y,\Fcal (Y))$ of type I with $H_4(Y,\CC)_\chi$ of dimension one.
According to  \cite{ls}, such a pair determines 
a $\G$-orbit $\Kcal_1$ of hyperplane sections of  $\PP(H_+)$.
Elsewhere on the GIT boundary,  the degeneration of the variation of polarized Hodge structure is  of  type II and so the period map is proper over the complement 
of $D(\Kcal_1)$ in  $\Gamma\bs \PP(H_+)$. The local Torelli theorem  says that $P$ is
a local isomorphism. Voisin's Torelli theorem \cite{voisin} says that if $X_1$ and $X_2$ are smooth cubic threefolds  with the same image under this period map, then their associated cubic fourfolds $Y_1,Y_2$ are projectively equivalent. We claim that for generic $X_1$, this isomorphism identifies the $\mu_3$-actions. The two actions certainly coincide on $H^{3,1}$ (it is there scalar multiplication), and hence they coincide on the transcendantal lattices of the two fourfolds  (the transcendental lattice of $Y_i$ is the smallest primitive sublattice in $H^4(Y_i)$ whose complex span contains 
$H^{3,1}(Y_i)$). But for generic $X_1$, the transcendental lattice is all of 
$H^4_\circ(Y_1)$. We conclude that $P$ is of degree one and hence is an open embedding. 

According to Lemma \ref{lemma:vanlattice}, some $K\in \Kcal_1$ contains a sublattice of $\Lambda$ isomorphic to $\Lambda_{10}$, in other words, $K\in \Hcal$. Since $\Kcal_1$ is a single $\G$-orbit,  it follows that $\Hcal=\Kcal_1$.

The image of $P$ is as asserted
if we prove that its image is disjoint with $D(\Hcal)$. If $P$ meets $D(\Hcal)$, then the complement of the image of $P$ is a closed subset of $\G\bs \PP(H_+)$ of codimension $\ge 2$ everywhere. Then $p$ will be an isomorphism and  hence identifies the GIT compactification with the Baily-Borel compactification. This is a contradiction since
the former has one dimensional boundary, whereas the latter's boundary is finite.
The isomorphism $P$ is covered by an isomorphism of line bundles: the $\G$-quotient of
$\Acal (3)|H_+^\circ$ gets identified with the orbifold  line bundle on 
$\GL (U)\bs S^\st$ defined by $\Ocal_{\PP (S)}(1)$ and this implies the last statement. 
\end{proof}

\begin{remarks}
The theory developed in \cite{looijenga} and \cite{ls} tells us a bit more. For example,
it interprets  $\proj \oplus_{k\ge 0} H^0(H^\circ_+, \Acal (k))^\G$ in terms of 
arithmetic data: the boundary it adds to $\G\bs \PP(H^\circ_+)$ has a  stratification 
whose members are described in Section \ref{sect:4folds}. Since this is also the GIT boundary,  which, according to Allcock, consists of an isolated point and a curve, we 
are able to recover these arithmetic data without any  extra effort: 
\begin{enumerate}
\item[(i)] $\G$ has two orbits of cusps in $\partial\PP(H_+)$,  corresponding to
the  cases ($D_4^3$) and ($A_5^2$) respectively.
\item[(ii)] A ($D_4^3$)-cusp does not lie on any  $\PP(K)$ with $K\in\Hcal$. 
\item[(iii)] The common intersection of the $K\in \Hcal$ for which $\PP(K)$ contains 
a given ($A_5^2$)-cusp is a codimension $2$  linear subspace of $H$.
\item[(iv)] The members of $\Hcal$ become disjoint when intersected with $H_+$.
\end{enumerate}
It  also tells us what the graph of the period map is like when regarded as a rational map from  the GIT compactification to the Baily-Borel compactification: on the GIT
side there is a hypersurface lying over the chordal cubic and on the Baily-Borel side 
there is a curve lying over the $(A_5^2)$-cusp. The hypersurface and the curve meet
in a single point.

If we are able to verify the first three of these four properties  in the context of lattice theory (which we have not done, though standard methods make this feasible), then the use of Voisin's Torelli theorem may be eliminated as follows: if we view the period map as a rational map from the GIT compactification to the Baily-Borel compactification, then (by the Zariski connectedness theorem)  the preimages of the two cusps will be disjoint. They are evidently contained in the GIT-boundary. One component of the GIT-boundary is the singleton represented by a cubic $3$-fold $X$ with three $D_4$-singularities. It gives rise to a  cubic $4$-fold $Y$ with three $\tilde E_6$-singularities and the period map sends this point to the ($D_4^3$)-cusp. Hence this singleton appears as a fiber of the rational period map.
It therefore suffices to verify that the latter is of degree one near $X$.
Now if $Y_t $ is a smooth cubic $4$-fold close to $Y$, then the orthogonal complement of the vanishing homology in $H^o_4(Y_t)$ is an isotropic plane and for this reason the behavior of the period map near $Y$ is essentially the one of its restriction to the vanishing homology, that is, the period map used in singularity theory
in the sense of \cite{looij:permap}. So we only  need to know that the latter is of degree one and this is indeed the case (see \emph{op.\ cit.}).
\end{remarks}

\section{The boundary pair defined by the chordal cubic}\label{sect:chordal}
In this section $X\subset \PP^4$ denotes a chordal cubic. That is, $X$ is the secant variety of a normal rational curve $C\subset \PP^4$ (since $C$ is unique up to projective equivalence, so is  $X$). This is indeed  a cubic hypersurface. A geometric argument might run as follows: if $\ell\subset \PP^4$ is a general line, then any point of $\ell\cap X$  lies on a secant  of $C$ by definition. If we project away from $\ell$ we map to a projective plane and the image of $C$ is an irreducible rational quartic curve in that plane. It will have three ordinary double points and these double points define the secants of $C$ which meet $\ell$. So  $\ell\cdot X=3$. 
Rather than making this argument rigorous, we derive an explicit equation for $X$ which is visibly of degree $3$: suppose $C$ be parameterized  by $[1:t]\mapsto [1:t:t^2:t^3:t^4]$.
Then for a general point $[x_0:\cdots :x_4]\in X$ there exist $s,t,\lambda,\mu$ such that 
$x_i=\lambda t^i+\mu s^i$. Elimination of $s,t,\lambda,\mu$ is straightforward and we find that 
\begin{equation}\label{eqn:chordal}
x_0(x_3^2-x_2x_4)+(x_ 2^3+x_1^2x_4-2x_1x_2x_3)=0
\end{equation}
is the cubic equation that defines $X$. Notice that $X$ contains the union $T_C$ of tangent lines of $C$ (as the secants of the infinitesimally near points).

\subsection*{Orbit decomposition of the chordal cubic}
It is convenient here to take a more abstract approach and construct everything in 
terms of a complex vector space $W$ of dimension two. Let us write $P_k$ for $\PP(\sym^kW)$, allowing  ourself to suppress the subscript when $k=1$.  We identify $P_k$ with the
 linear system of effective degree $k$ divisors on $P$ and  identify $P$  with its image  in $P_k$ (for $k>0$) by means of $p\in P\to kp\in P_k$.

Our curve $C$ is now  identified with $P\subset P_4$ and $X\subset P_4$ parameterizes the quartics in $P$ that can be written as the sum of two fourth powers, or at least  infinitesimally so (these are quartic forms that are divisible by a third power). In the last case this means that the point in question lies on a tangent line of $C$, i.e., in $T_C$.
If we interpret $P_4$ as the linear system of effective degree four divisors $x$ on 
$P$, then $x\in X$ if and only if
it is invariant under an involution which fixes at least one point of $\supp (x)$.  
We then see that $\PGL(W)$ has three orbits in $X$: $C$ (divisors of the form $4p$), $T_C-C$ (divisors of the form $3p+q$ with $q\not=p$) and $X-T_C$ (reduced divisors).
The singular locus of $X$ is $C$.

\subsection*{Cohomology of the smooth part}
We use the classical theory of Lefschetz  pencils to determine the  $\chi$-Betti numbers of $Y_\reg$.  Both $X$ and $Y$ have a smooth singular locus with a transversal singularity whose link is a rational homology sphere. So either is a rational homology manifold and therefore its cohomology satisfies Poincar\'e duality over $\QQ$.
The Gysin sequence for the pair  $(Y,C)$,
\[
\cdots\to H^{k-4}(C;\QQ)\to H^k(Y;\QQ)\to H^k(Y_\reg;\QQ)\to H^{k-3}(C;\QQ)\to\cdots ,
\]
shows that $H^k(Y;\CC)_\chi\to H^k(Y_\reg;\CC)_\chi$ is an isomorphism (and likewise for $\bar\chi$, of course), so that either comes with a nondegenerate $\QQ(\zeta)$-valued Hermitian form.
Let us do now an Euler characteristic computation. 
The projectively completed tangent bundle $T_C$ of $C$ has Euler characteristic equal to $e(C).e(\PP^1)=4$. Since $X-T_C$ is a fibered (over $P_2-P$) with fiber $\CC^\times$, its Euler characteristic is zero. Hence $e(X)=4$ and $e(P_4-X)= 5-4=1$. Since $Y\to P_4$ is a $\mu_3$-cover totally ramified over $X$, we have $e(Y)=3e(P_4-X)+e(X)=3+4=7$.
The $\mu_3$-orbit space of $Y$ can be identified with $P_4$, which has Euler characteristic $5$. Hence $e_\chi(Y)=e_{\bar\chi}(Y)=\half (7-5)=1$.

In order to prove the  more precise Proposition \ref{prop:chordal} we choose a Lefschetz pencil for $X$. Such a pencil is defined  a generic $2$-plane
$A\subset P_4$ (its axis). The genericity assumption  entails that
it  avoids $C$ and meets $X_\reg$  transversally. A generic hyperplane $H$ through $A$ meets $X$ in a cubic surface 
$X_H$ whose singular points are $C\cap H$. The latter intersection is transversal and so
each of the four points of $C\cap H$ is a singularity of type $A_1$. 
The  cubic surfaces with four $A_1$-singularities form a single projective
equivalence class. 

\begin{lemma}\label{lemma:chicomp1}
Let $Z\subset \PP^3$ be a cubic surface with $4$ $A_1$-singularities and let
$K\to \PP^3$ be the normal $\mu_3$-cover that totally ramifies over $Z$. Then
$H^\pt (K;\CC)_\chi =H^3 (K;\CC)_\chi$  and the latter has dimension one.
\end{lemma}
\begin{proof}
We first observe that $K$ is  a cubic $3$-fold with $4$ $A_2$-singularities. Since the vanishing homology of an $A_2$-singularity is a symplectic  unimodular lattice of rank 
$2$, we find that  the primitive cohomology of $K$ is concentrated in degree $3$ and is of rank $4.2=8$ less than the rank of the primitive cohomology of a smooth cubic in that dimension (which is $10$). It follows that  $H_\pt (K;\CC)_\chi=H_3 (K;\CC)_\chi$ has dimension one.
\end{proof}

\begin{lemma}\label{lemma:tgt}
If $H$ is tangent hyperplane  of $X$ at some $p\in X-T_C$, then $H^\pt (Y_H;\CC)_\chi$ is trivial.
\end{lemma}
\begin{proof}
Since $\SL (W)$ is transitive on $X-T_C$, any $p\in X-T_C$ will do. We use the equation  \ref{eqn:chordal} of $X$, $x_0x_3^2-x_0x_2x_4+ x_2^3+x_1^2x_4-2x_1x_2x_3$, and we take $p=[1:0:0:0:1]$.
Then the tangent hyperplane at $p$ is given by $x_2=0$ and so $X_H$ is given by
$x_0x_3^2+x_1^2x_4$. Hence $Y_H$ is given by $y^3=x_0x_3^2+x_1^2x_4$. 
The singular set of this $3$-fold is the line $\ell$ defined by $y=x_1=x_3=0$. 
Consider the $\CC^\times$-action on $Y_H$ defined by 
\[
\lambda [x_0:x_1:x_3:x_4:y]:= [x_0:\lambda^3 x_1:\lambda^3 x_3: x_4:\lambda^2 y]
\]
We notice that the fixed point set is the union of $\ell$ and the line $\ell'$ defined by
$y=x_0=x_4=0$. This action provides a contraction of $Y_H-\ell'$ onto $\ell$ so that
$\ell\subset Y_H-\ell'$ is a homotopy equivalence. This shows that 
$H^\pt(Y_H-\ell';\CC)_\chi=0$. We  also have $H^\pt(\ell';\CC)_\chi=0$, of course.
Since $\ell'$ lies in the smooth part of $Y_H$, the Gysin sequence can be applied to the pair  $(Y_H,\ell')$. We thus find that $H^\pt(Y_H;\CC)_\chi=0$. 
\end{proof}

\begin{proof}[Proof of Proposition \ref{prop:chordal}]
Since $Y$ is a rational homology manifold, $H^\pt(Y;\CC)_\chi$ satisfies Poincar\'e duality in the sense that $H^\pt(Y;\CC)_\chi$ pairs nondegenerately with 
$H^{8-\pt}(Y;\CC)_{\bar\chi}$. In particular, $H^4(Y;\CC)_\chi$ comes with a nondegenerate Hermitian form. As to the assertions concerning $H^\pt(Y:\CC)_\chi$,  
since we already verified that $e_\chi (Y)=1$, it suffices to prove that
$H^k(Y;\CC)_\chi=0$ for $k\le 3$.

Let $L_A$ be the pencil of hyperplanes in $\PP^4$ passing through $A$ and  
denote by $\tilde X\subset \PP^4\times L_A$ resp.\ $\tilde Y\subset \PP^5\times L_A$ the 
corresponding Lefschetz pencils. The projection $\tilde Y\to Y$ contracts  
$Y_A\times L_A$ along its projection onto $L_A$ and so 
$H^\pt(Y;\CC)\to H^\pt(\tilde Y;\CC)$ is injective. We prove that
$H^k(\tilde Y;\CC)_\chi=0$ for $k\le 3$. To this end, consider the $\chi$-Leray spectral sequence
of the projection $\pi :\tilde Y\to L_A$. Lemma \ref{lemma:chicomp1} shows that 
$(R^q\pi_*\CC_{\tilde Y})_\chi$ is zero unless $q=3$ so that
the Leray spectral sequence for $\chi$-cohomology degenerates and 
\[
H^k(\tilde Y;\CC)_\chi=H^{k-3}(L_A;(R^3\pi_*\CC_{\tilde Y})_\chi). 
\]
In particular $H^k(\tilde Y;\CC)_\chi=0$ for $k\le 2$.
Lemma \ref{lemma:tgt} implies that  a stalk of $(R^q\pi_*\CC_{\tilde Y})_\chi$
is zero if the associated hyperplane is tangent to $X_\reg$. Since there 
are such hyperplanes, 
$H^3(\tilde Y;\CC)_\chi=H^0(L_A;(R^3\pi_*\CC_{\tilde Y})_\chi)=0$  also.
The rest of the proposition follows from an application of Lemma 
\ref{lemma:type1}.
\end{proof}

\end{document}